\documentstyle[12pt,twoside]{amsart}

\title
{The Local Nash problem on arc families of singularities}
\author{Shihoko Ishii} 
\address{Department of Mathematics, Tokyo Institute of
Technology, Oh-Okayama, Meguro, 152-8551, Tokyo, Japan
\newline
e-mail : shihoko@@math.titech.ac.jp}

\newcommand{\bC}{{\Bbb C}}

\newcommand{\bZ}{{\Bbb Z}}

\newcommand{\bR}{{\Bbb R}}

\newcommand{\bm}{{\bf m}}

\newcommand{\Spec}{\operatorname{Spec}}

\newcommand{\Hom}{\operatorname{Hom}}
\newcommand{\Sing}{\operatorname{Sing}}
\newcommand{\ord}{\operatorname{ord}}

\newcommand{\st}{{\Spec k[[t]]}}

\newcommand{\sT}{{\Spec K[[t]]}}

\renewcommand{\o}[0]{{\mathcal O}} 

\newcommand{\spec}[0]{\operatorname{Spec}}
\newcommand{\sing}[0]{\operatorname{Sing}}

\def\to {\longrightarrow}

\newtheorem{thm}{Theorem}[section]

\newtheorem{lem}[thm]{Lemma}
\newtheorem{cor}[thm]{Corollary}
\newtheorem{prop}[thm]{Proposition}

\theoremstyle{definition}
\newtheorem{defn}[thm]{Definition}

\newtheorem{say}[thm]{}
\newtheorem{exmp}[thm]{Example}

\newtheorem{rem}[thm]{Remark}

\theoremstyle{remark}

\begin{document}
\maketitle
\footnote{partially supported by Grand-In-Aid of Ministry of Science 
and Education in Japan}

\begin{abstract}
  This paper shows the affirmative answer to the local Nash problem
  for a toric singularity and analytically pretoric singularity.
  As a corollary we obtain the affirmative answer to the local 
  Nash problem for 
  a  
  quasi-ordinary singularity. 
\end{abstract}

\section{Introduction}
\noindent
  The Nash problem was posed by John F. Nash in his preprint in 1968, 
  which was later published as \cite{nash} in 1995.  
  The problem in his paper reads in two ways:
  
  \begin{enumerate}
  \item
  the bijectivity of the map from the set of the families of arcs 
  passing through  ``the singular locus'' to the set of the essential 
  divisors 
  over ``the singular locus''
  \item
  the bijectivity of the map from the set of the families of arcs 
  passing through ``a singular point'' to the essential divisors over 
  ``the singular point''
\end{enumerate}  
  For convenience sake, we call the former the Nash problem and the latter 
  the local Nash problem.
  For a variety with an isolated singularity, the two problems coincide.
  
  In case of a 2-dimensional normal (therefore isolated) 
  singularity, the problem is studied in \cite{LR}, \cite{P},  \cite{RL}
   The Nash problem for general dimension is studied in \cite{i2}, \cite{i-k}.
   (More detailed information about the known facts will be given in 
   \ref{knownfacts} in the second section.)
  In this paper we study the local Nash problem.
  We show the affirmative answer to the local Nash problem for 
  every point of a toric variety and also for an analytically pretoric 
  singularity.
  As a corollary we obtain the affirmative answer to the local Nash 
  problem for a 
   quasi-ordinary singularity.
  
  This paper is organized as follows:
  In the second section, we introduce the Nash map, the local Nash 
  map, the Nash 
  problem and the local Nash problem.
  In the third section, we show the affirmative answer to the local 
  Nash problem for a toric variety.
  In the fourth section we show the affirmative answer to the local 
  Nash problem for an analytically pretoric  singularity.
  As a corollary we obtain the affirmative answer for a 
   quasi-ordinary singularity.

  In this paper we work on schemes  over an algebraically closed field \( 
  k  \) of arbitrary characteristic.
  All \( k \)-schemes are assumed to be  pure dimensional excellent 
  schemes over \( k \).
  All reduced \( k \)-schemes are moreover assumed to have open dense 
  regular locus.
  By a regular \( k \)-scheme we mean a  \( k \)-scheme with every
  local ring regular.
 

\section{The Nash problem and the local Nash problem}

\noindent
\begin{defn}
  Let \( X \) be a scheme  over \(k\)
and $K\supset k$ a field extension.
  A morphism  \(\alpha: \Spec K[[t]]\to X \) is called an {\it arc} of \( X \).
  We denote the closed point of \( \Spec K[[t]] \) by \( 0 \) and
  the generic point by \( \eta \).
\end{defn}

  For a \( k \)-scheme \( X \),
  the arc space \( X_{\infty} \)
  is characterized by the 
  following property (\cite{v}):  
 \begin{prop}
\label{arcspace}
  Let \( X \) be a \( k \)-scheme.
  Then
  \[ \Hom _{k}(Y, X_{\infty})\simeq\Hom _{k}(Y\widehat\times_{\Spec k}
  \st, X) \]
   for an arbitrary \( k\)-scheme \( Y \),
   where \( Y\widehat\times_{\Spec k}\st \) means the formal completion 
   of \( Y\times_{\Spec k}\st \) along the subscheme 
   \( Y\times _{\Spec k} \{0\} \).
\end{prop}

\begin{say}
  By thinking of the case \( Y=\Spec K \) for an extension field \( K \) 
  of \( k\), 
  we see that  \( K \)-valued points of \( X_{\infty} \) correspond to
  arcs \( \alpha: \sT\to X \) bijectively. 
  Based on this, we denote the \( K \)-valued point corresponding to 
  an arc \( \alpha: \sT\to X \) by the same symbol \( \alpha \).
    The canonical projection \( X_{\infty} \to X \), \( \alpha\mapsto 
  \alpha(0) \) is  denoted by \( \pi_{X} \). 
 
  A morphism \( \varphi:Y \to X \) of varieties induces a canonical 
  morphism \( \varphi_{\infty}:Y_{\infty} \to X_{\infty} \), \( 
  \alpha \mapsto \varphi\circ\alpha \). 
\end{say}

\begin{defn}
\label{nash}
  Let \( X \) be a reduced \( k \)-scheme and \( \Sing X \) 
  the singular locus of \( 
  X \), i.e., the set of the points whose local rings are not regular.
  Recall that  we assume that all reduced \( k \)-schemes are pure 
  dimensional  excellent schemes and have the open dense regular locus. 
  An irreducible component  \( C \) of \( \pi_{X}^{-1}(\Sing X) \) is called a 
 {\it Nash component} of \( X \) if \( C  \) is not contained in \( (\Sing X)_{\infty} 
  \).
   (In \cite{i-k} a Nash component is called a ``good component''.) 
  Let \( x \) be a (not necessarily closed) point of  \( X \).
  An irreducible component \( C \) of  \( \overline{\pi_{X}^{-1}(x)} \)
   is called a 
  {\it local Nash component} of \( (X,x) \) if \( C \) is not 
  contained in \( (\Sing X)_{\infty} 
  \).
\end{defn}
  Here, we note that  
   every irreducible component of \( \pi_{X}^{-1}(\Sing 
  X) \) 
  is a Nash component if \( k \) is of characteristic zero
  (\cite[Lemma 2.12]{i-k}).

\begin{say}
\label{key}
  Assume that  \( \varphi: Y \to X \) is a proper morphism which is an  isomorphism away from
  \( \Sing X \).
  Let \( \alpha \) be the generic point of 
   a Nash component or of a local Nash component.
  Then \( \alpha(\eta)  \) is outside of \( \Sing X \), 
  therefore it is lifted to \( Y \) by the isomorphism \( \varphi \).
  Then, by the valuative criterion of properness 
  \( \alpha \)
  can be uniquely lifted to an arc of \( Y \).
  This property is essential for our arguments in this paper.
\end{say} 
  
\begin{lem}
\label{fiber}
  Let \( X \) be an integral \( k \)-scheme and \( x \) an 
  analytically irreducible point of \( X \), 
  i.e., \( \widehat{\cal O}_{X,x} \) 
  is an integral domain.
  Let \( \widehat{X} \) be \( \spec \widehat{\cal O}_{X,x} \).
  Then, the canonical morphism \( \iota_{\infty}: \widehat{X}_{\infty}
  \to 
  X_{\infty} \) induces an isomorphism 
   \( \pi_{\widehat{X}}^{-1}(x)\simeq \pi_{{X}}^{-1}(x) \),
  where the closed point of \( \widehat{X} \) is also denoted by \( x \).
\end{lem}  

\begin{pf}
  First, note that the canonical morphism \( \iota:\widehat{X}\to
  X \) gives the morphism \( \iota_{\infty}:\widehat{X}_{\infty}\to 
  X_{\infty} \) whose restriction gives 
  \( \iota_{\infty}:\pi_{\widehat{X}}^{-1}(x)\to \pi_{X}^{-1}(x) \).
  We may assume that \( X=\spec A \) for a \( k \)-algebra \( A \).
 Let \( \pi_{X}^{-1}(x)=\Spec R \). 
  By Proposition \ref{arcspace},
  the inclusion  \( \pi_{X}^{-1}(x)\subset X_{\infty} \) induces a
  homomorphism \( {\cal O}_{X,x} \to R[[t]]\)
  which sends the maximal ideal of \( {\cal O}_{X,x} \) to the ideal 
  \( (t) \).
  Then, we get the homomorphism of projective limits 
  \[ \begin{array}{ccc}
   \widehat{\cal O}_{X,x} &  \to  &  R[[t]] \\
   \| & &  \|  \\
   \varprojlim\o_{X,x}/{\frak{m}}_{X,x}^m  & &  
   \varprojlim 
  R[[t]]/(t)^m .\\
  \end{array} 
  \]
  
   Again by Proposition  \ref{arcspace},
  this homomorphism gives a morphism 
  \( \pi_{X}^{-1}(x)\to \widehat{X}_{\infty} \) whose image is in 
  \( \pi_{\widehat{X}}^{-1}(x) \).
  This is the inverse morphism of \( 
  \iota_{\infty}:\pi_{\widehat{X}}^{-1}(x)\to
  \pi_{X}^{-1}(x) 
   \).    
\end{pf}  
  
\begin{lem}
\label{irred}
  Let \( X \) be a regular \( k \)-scheme and \( E \) an irreducible 
  regular
  closed subset of \( X \).
  Then \( \pi_{X}^{-1}(E) \) is an  irreducible  closed subset of 
  \( X_{\infty} \).
  \end{lem}  
  
\begin{pf}
  We may assume that \( X=\spec A \) for an integral domain \( A \).
  As \( {\cal O}_{X,p} \) is a regular local ring for every \( p\in X \),
  we have   \( \widehat{{\cal O}}_{X, p}=k(p)[[x_{1},\ldots,x_{n}]] \)
  for  some indeterminates \( x_{1},\ldots,x_{n} \),
  where \( k(p) \) is the residue field of \( \o_{X,p} \).
  If we put \( \widehat{X}=\spec \widehat{{\cal O}}_{X, p}\),
  this shows that 
   \( \pi_{X}^{-1}(p)=\pi_{\widehat{X}}^{-1}(p) 
  \) is irreducible for every \( p\in X \).
  Therefore, it is sufficient to prove that 
  \( \pi_{X}^{-1}(p)\subset 
  \overline{\pi_{X}^{-1}(q)} \)
  for \( p,q \in X \) with \( p\in \overline{\{q\}} \) and 
   \( \overline{\{q\}} \) regular.
  Let \( \frak{p} \) and \( \frak{q} \) be the prime ideals in \( 
  \o_{X,p} \)
  corresponding to \( p \) and \( q \), respectively.
  Then, we may assume that \( {\frak{p}}=(x_{1},..,x_{r},x_{r+1},..,x_{n}) \) and 
  \( {\frak{q}} =(x_{1},..,x_{r})\).
  Let \( \alpha \) be the generic point of \( \pi_{X}^{-1}(p) \).
  Then \( \alpha \) induces a local homomorphism 
  \( {\cal O}_{X,p}\to K[[t]] \) and this can be extended to a local 
  homomorphism 
 \[   \alpha^*:\widehat{{\cal 
  O}}_{X,p}=k(p)[[x_{1},..,x_{r},x_{r+1},..,x_{n}]]\to K[[t]] .  
 \]
 Define  \[ \Lambda^*:k(p)[[x_{1},..,x_{r},x_{r+1},..,x_{n}]]\to 
 K[[\lambda_{r+1},\ldots,\lambda_{n}, t]]  \] 
 by 
 
 \begin{enumerate}
 \item[ ]
 \( \Lambda^*(x_{i})= \alpha^*(x_{i}) \) for 
 \( i=1,..,r \) and 
\item[ ] 
 \( \Lambda^*(x_{i})=\lambda_{i}+ \alpha^*(x_{i}) \)
 for \( i=r+1,..,n \).
\end{enumerate} 
  Here \( \lambda_{r+1},\ldots,\lambda_{n} \) are indeterminates. 
  The restriction of this map onto \( A \) gives a family of arcs 
  \( \Lambda:\spec K[[\lambda_{r+1}, \ldots, \lambda_{n}]] \to  X_{\infty} \).
  Let \( 0' \) and \( \eta' \) be the closed point and the generic 
  point of   \( \spec K[[\lambda_{r+1},\ldots,\lambda_{n}]] \), respectively.
  Denote the quotient field of \( K[[\lambda_{r+1}, \ldots, \lambda_{n}]] \)
  by \( K((\lambda_{r+1},  \ldots,  \lambda_{n})) \).
  Then \( \Lambda(0')=\alpha \) and \( \beta:=\Lambda(\eta'):
  \spec K((\lambda_{r+1},\ldots,\lambda_{n}))[[t]]\to X \) is an arc in \( \pi_{X}^{-1}(q) \),
  since \( {\beta^{*-1}}((t))={\frak{q}}\cap A \).
  This yields that \( \alpha\in \overline{\{\beta\}}\subset 
  \overline{ \pi_{X}^{-1}(q)} \).    
\end{pf}

We note that if \( X \) is a non-singular variety, 
\( \pi_{X}^{-1}(E) \) is always irreducible for an irreducible subset
  \( E \).

\begin{defn}
  A birational morphism \( \varphi:Y \to X \) of reduced \( k 
  \)-schemes is a morphism which gives a bijection between the sets 
  of the irreducible components of \( Y \) and \( X \) and the 
  restriction of \( \varphi \) on each irreducible component is birational.

 Let $X$ be a reduced \( k \)-scheme, 
\( \psi:X_1 \to X \)
a proper birational morphism from a normal  \( k \)-scheme 
\( X_{1} \)
 and $E\subset X_1$ 
an irreducible exceptional divisor  of 
  \( \psi \). 
 Let  \( \varphi:X_2 \to X \) be another
 proper birational morphism from a normal  \( k \)-scheme 
 \( X_{2} \).
The birational  map
  \(\varphi^{-1}\circ \psi  : X_1  \dasharrow  X_2 \) is 
defined on a (nonempty) open subset $E^0$ of $E$.
The closure of $(\varphi^{-1}\circ \psi)(E^0)$ is 
  called the {\it center} of $E$ on $X_2$.

  We say that \( E \) appears in \( \varphi \) (or in \( X_2 \)),
  if 
the center of $E$ on $X_2$ is also a divisor. In this case
the birational  map
  \(\varphi^{-1}\circ \psi  : X_1  \dasharrow  X_2 \) is  a local isomorphism at the 
  generic point of \( E \) and  we denote the birational
  transform of \( E \) on \( X_2 \)
  again by \( E \). For our purposes $E\subset X_1$ is identified
with $E\subset X_2$. 
 Such an equivalence class  is called an {\it exceptional divisor over $X$}. 
  
  An exceptional divisor \( E \) over \( X \) is called an 
  {\it exceptional divisor over} \( (X,x) \) for a point \( x\in X \) if 
  the center of \( E \) on \( X \) is \( \overline{\{x\}} \).

\end{defn}

\begin{defn}
  Let \( X \) be a reduced \( k \)-scheme. 
  In this paper, by  a
{\it  resolution} of the 
  singularities of \( X \) 
   we mean  a proper birational  morphism \( \varphi:Y\to X \) 
  with   a regular  \( k \)-scheme \( Y \) 
  such that 
  the restriction $Y\setminus \varphi^{-1}(\sing X)\to X\setminus \sing X$
 is an isomorphism.
 
 The existence of a resolution for  a reduced \( k 
 \)-scheme \( X \) is a difficult problem.
 For a variety over a field of characteristic zero 
 the existence of a resolution was proved by Hironaka \cite{hironaka}.
 But for a general  reduced \( k \)-scheme it is still an open 
 problem.
 From now on, we always assume the existence of a resolution.

 \end{defn}
  
  \begin{defn}\label{nashessential}
  An exceptional divisor \( E \) over a reduced 
  \( k \)-scheme \( X \) is called an {\it essential 
  divisor} over \( X \)
   if for every resolution \( \varphi:Y\to X \)
  the center of \( E \) on \( Y \) is an irreducible component of 
   \( \varphi^{-1}(\sing X) \).
   The center of an essential divisor over \( X \) on a resolution \( 
   Y \) is called an {\it essential component on } \( Y \).
  
  For  a point \( x\in 
  X \)  an exceptional divisor \( E \) over \( (X,x) \) is called 
  an {\it essential divisor over} \( (X,x) \) if for every resolution
  \( \varphi:Y\to X \)
   the center of \( E \) on \( Y \) is an irreducible component of 
  \( \overline{\varphi^{-1}(x)} \).
  The center of an essential divisor over \( (X,x) \) on a resolution \( Y \) 
  is called an {\it essential component over } \( (X,x) \) on \( Y \).
\end{defn}

\begin{rem}
  Take an integral scheme \( X \) and a point \( x\in X \).
  There are canonical bijections: 
  \[ \{\mbox{essential divisors over \( X \)}\}\simeq
  \{\mbox{essential components on a resolution \( Y \)}\},  \]
  \(\{\mbox{essential divisors over \( (X,x) \)}\}\)
  
  \hskip2truecm
  \(\simeq
  \{\mbox{essential components over \( (X,x) \) on a resolution \( Y 
  \)}\}. \)
  Indeed, for an essential divisor \( E \), let \( \Phi(E) \) be the 
  center of \( E \) on \( Y \).
  Then we have a map \( \Phi \) from the set of essential divisors to the set of 
  the essential components.
  Conversely, for an essential component \( C \) on \( Y \),
  take the blow-up \( \tilde Y\to Y \) with the center \( C \) 
  and let \( E \) be the unique exceptional divisor which is mapped onto \( C \).
  Then \( E \) is an essential divisor whose center on \( Y \) is \( C \).
\end{rem}

\begin{say}
{\bf The Nash problem}
 Let \( \varphi:Y \to X \) be a resolution of the singularities of a 
 reduced 
   \( k \)-scheme \( X \) such that  \( \varphi^{-1}(\Sing X) \) 
  is a union of non-singular divisors. 
  Let \( \varphi^{-1}(\Sing X)=\bigcup_{j} E_{j} \) be the decomposition 
  into irreducible components.
  Let \( \{C_{i}\} \) be the Nash components of \( X \).
  Then the morphism \( \varphi_{\infty}: \bigcup_{j} \pi_{Y}^{-1}(E_{j}) \to 
  \bigcup_{i}C_{i} \) is dominant and bijective outside  
  \( (\Sing X)_{\infty} \) by \ref{key}.
  As \( \pi_{Y}^{-1}(E_{j}) \)'s are irreducible by 
  Lemma \ref{irred}, for each \( C_{i} \) there 
  is unique \( E_{j_{i}} \) such that 
  \( \pi_{Y}^{-1}(E_{j_{i}}) \) is dominant over \( C_{i} \).
  In \cite{nash} Nash proved that this  \( E_{j_{i}} \) is an essential divisor 
  over \( X \) 
  (for the proof see also \cite[Theorem 2.15]{i-k}).
  This map
  \begin{center}
   \( {\cal N}: \{ \) Nash components \( \}\to \{ \) essential 
  divisors over \( X \}\), \( C_{i } \mapsto E_{j_{i}} \) 
  \end{center}
  is called {\it the Nash map}.
  Obviously this map is injective and the Nash problem asks if 
  this map is bijective.
\end{say}

\begin{say}
\label{locNash}
{\bf The local Nash problem}
Let \( \varphi:Y \to X \) be a resolution of the singularities of a 
 reduced \( k \)-scheme such that  
\( \overline{\varphi^{-1}(x)} \)  
  is a union of non-singular  divisors. 
  Let \( \overline{\varphi^{-1}(x)}=\cup_{j}E_{j} \) be the 
  decomposition into irreducible components. 
   Let \( \{C_{i}\} \) be the local Nash components of \( (X,x) \).
    Then the morphism \( \varphi_{\infty}: \bigcup_{j} \pi_{Y}^{-1}(E_{j}) \to 
  \bigcup_{i}C_{i} \) is dominant and  injective outside \( (\Sing 
  X)_{\infty} \) by \ref{key}.
  As \( \pi_{Y}^{-1}(E_{j}) \)'s are irreducible, for each \( C_{i} \) 
  there is a unique \( E_{j_{i}} \) such that 
  \( \pi_{Y}^{-1}(E_{j_{i}}) \) is dominant to \( C_{i} \).
  By the following lemma, this  \( E_{j_{i}} \) is an essential divisor 
  over \( (X,x) \).
   This map
  \begin{enumerate}
  \item[ ]
   \({{\ell} {\cal N}}:\{ \) local Nash components of \((X,x) \}\)
   \item[ ]
   \ \ \ \ \ \ \ \ \ \(\to \{ \) essential 
  divisors over \(( X , x) \}\), \( C_{i } \mapsto E_{j_{i}} \) 
  \end{enumerate}
  is called {\it the local Nash map}.
  Obviously this map is injective and the local Nash problem asks if 
  this map is bijective.
  
  If \( x\in X \) is a unique singularity on \( X \),
  then the Nash problem for \( X \) is the same as the local Nash 
  problem for \( (X,x) \).
\end{say}

\begin{lem}
 Under the notation above, \( E_{j_{i}} \) is an essential divisor over 
 \( (X,x) \).
\end{lem}

\begin{pf}
  Let \( \psi:Y'\to X \) be any resolution.
  Let \( E_{j_{i}}' \) be the center of \( E_{j_{i}} \) on \( Y' \).
   Then, \( C_{i}
  =\overline{\psi_{\infty}\pi_{Y'}^{-1}(E_{j_{i}}')} \).
  Let \( E' \) be an irreducible component  of \( 
  \overline{\psi^{-1}(x)} \) containing \( E_{j_{i}}' \).
  Then
  \[ C_{i} =\overline{\psi_{\infty}\pi_{Y'}^{-1}(E_{j_{i}}')}
  \subset \overline{\psi_{\infty}\pi_{Y'}^{-1}(E')  },
  \] 
  where the last term is in \( \overline{ \pi_{X}^{-1}(x) } \).
  As \( C_{i} \) is an irreducible component of 
  \( \overline{ \pi_{X}^{-1}(x) } \),
  the above inclusion is the equality.
  By the bijectivity of \( \psi_{\infty} \) outside \( (\sing 
  X)_{\infty} \) 
  the generic points \( \alpha \) and \( \alpha' \) of \( \psi_{\infty}\pi_{Y'}^{-1}(E_{j_{i}}') \) and 
  \( \psi_{\infty}\pi_{Y'}^{-1}(E') \), respectively, must coincide,
  which yields that the generic points of \( E_{j_{i}}' \) and 
  \( E' \) coincide,
  because \( E_{j_{i}}'=\overline{\{\alpha(0)}\} \) and 
  \( E'= \overline{\{\alpha'(0)}\} \).
\end{pf}

\begin{say}
\label{knownfacts}
{\bf Known facts on the Nash problem.}
   An essential divisor, which is a slightly different notion from ours, is studied by Catherine Bouvier and G\'erard 
   Gonzalez-Sprinberg in \cite{BG}.
  The idea of the proof of a theorem in this paper is very useful for our discussion.
  The Nash problem is affirmatively answered for \( A_{n} 
  \)-singularities by John F. Nash \cite{nash}, for a minimal 
  singularity on a surface  by Ana Reguera \cite{RL} and for a sandwiched surface 
  singularity By Monique Lejeune-Jalabert and Ana Reguera \cite{LR},  
  \cite {RL2}.
  Recently the author was announced that the affirmative answer is 
  proved for a 
  \( D_{n} \)-singularity on a surface by Camille Plenat.
  Camille Plenat and Popescu-Pampu \cite{P} proved the affirmative answer to
  certain non-rational singularities with  combinatorial conditions.
  The Nash problem is  affirmatively answered also for a toric variety of 
  arbitrary dimension  in
  \cite{i-k}. 
  But affirmative answer does not hold for a general singularity.
  The same paper \cite{i-k} gives a counter example of dimension 4,
  therefore we have counter examples for dimension higher than 4
  by making the product with a non-singular variety.
  For dimension 2 and 3 the problem is still open.
  These are all for a normal variety.
  We should note that, this problem for a non-normal variety 
  is not automatically reduced to the case of  
  the normalized variety. 
  In spite of that, for a non-normal toric variety the Nash problem is 
  affirmatively proved in (\cite{i2}).
  A non-normal toric variety has much stronger properties than
  just the fact that its normalization is a toric variety.

\noindent
{\bf Known facts on the local Nash problem.}
  As a normal surface singularity is isolated,
  all results on the Nash problem for a normal surface singularity 
  are the results on the local Nash problem.
  The counter example to the Nash problem given in \cite{i-k} is 
  an isolated singularity, 
  therefore it is also a counter example to the local Nash problem.
  Hence, the next step to study is to know in which category the local 
  Nash problem (or the Nash problem) is affirmative.
\end{say}

Now we close this section with the following basic lemma, 
which implies that a Nash component and a local Nash component are 
``fat''  in terms of \cite{i2}.

\begin{lem}
\label{fat}
  Let \( C \) be a Nash component of an integral \( k \)-scheme \( X \) or a local 
  Nash component of \( (X,x) \) for a point \( x \) of 
  an integral \( k \)-scheme \( X \).
  Let \( \alpha :\sT \to X \) be the generic point of \( C \).
  Then, \( \alpha(\eta) \) is the generic point of \( X \), which is 
  equivalent to that the corresponding ring homomorphism 
  \( \alpha^*:\Gamma(U, \o_{X})\to K[[t]] \) is injective, 
  where \( U \) is an 
  affine open neighborhood of \( \alpha(0) \). 
\end{lem} 

\begin{pf}
  We prove the statement for a local Nash component.
  The other case is essentially the same.
  Let \( C \) be a local Nash component of \( (X,x) \)
  and \( E \)  an  essential divisor over \( (X,x) \) corresponding to \( C \).
  Let \( \varphi:Y \to X \) be a resolution of the singularities of 
  \( X \), on which the divisor \( E \) appears.
  Then \( \overline{\varphi_{\infty}(\pi_{Y}^{-1}(E) )}=C \).
  As \( \pi_{Y}^{-1}(E) \) is an irreducible cylinder on a 
  non-singular variety \( Y \), it is not contained in the arc space 
  of any proper closed subscheme of \( Y \).
  Therefore, the generic point \( \beta \) of \( \pi_{Y}^{-1}(E) \)
  sends the generic point of \( \sT \) to the generic point of 
  \( Y \).
  Hence, the generic point \( \varphi_{\infty}\beta \) of \( C \) also
  sends the generic point of \( \sT \) to the generic point of 
  \( X \).
\end{pf}

\section{The local Nash problem for a toric variety}

\noindent
  In this section we prove the local Nash problem for a toric variety.
  First we remark some basic notion of the arc space of a 
  toric variety.
   Here, we use the notation and terminologies of \cite{fulton}.
  Let  $M$  be the free abelian group  ${\bZ}^n$ $(n\geq 1)$
  and  $N$   its dual $\Hom_{\bZ}(M, {\bZ})$.
  We denote  $M\otimes _{\bZ}{\bR}$  and $N\otimes_{\bZ}{\bR}$  by
  ${M_{\bR}}$  and  $N_{\bR}$, respectively.
  The canonical pairing \( \langle\ , \ \rangle:
  N\times M \to \bZ \)  extends to  
  \( \langle\ \ , \ \rangle:
  N_{\bR}\times M_{\bR} \to \bR \).
   The group ring \( \bC[M] \) is generated by monomials \( x^{\bm} \) 
  (\( \bm \in M \)) over \( \bC \).
  A {\it cone in} \( N \) is a strongly convex rational polyhedral cone in \( N_{\bR} \).
  
\begin{say}
\label{T} 
  Let \( X \) be an affine toric variety defined by a cone \( \sigma \) 
  in \( N \).
  In \cite{i}, for \( v\in \sigma\cap N \) we define 
 \[  T^X_{\infty}(v)= \{\alpha\in X_{\infty} \mid \alpha(\eta)\in T,\  \ord_{t} 
  \alpha^*(x^u)=\langle v, u \rangle \ {\operatorname{for}}\ u\in M \}, \]
   where \( T \) denotes the open orbit and also the torus acting on \( X \).
  The set \(  T^X_{\infty}(v) \) is an irreducible locally closed 
   subset of \( X_{\infty} \) which is not contained in \( (\sing 
   X)_{\infty} \)(\cite{i}).

   Let \( \tau \) be the face of \( \sigma \) such that 
   \( v\in \tau^o \), where \( \tau^o \) means the relative interior of 
   \( \tau \).
   Then for every \( \alpha\in T_{\infty}^X(v) \), we have that \( \alpha(0)\in 
   orb (\tau) \) (\cite[Proposition 3.9]{i-k}).
\end{say}

  In \( \sigma\cap N \) we define an order \( \leq_{\sigma} \) as 
  follows:
  \[ v\leq_{\sigma}v' \Leftrightarrow v'-v\in \sigma. \]
  Then the following is obtained in \cite{i}.
   
\begin{prop}[\cite{i}]
\label{dominate}
  Let \( X \) be an affine toric variety defined by a cone \( \sigma \) 
  in \( N \).
  For \( v,v' \in \sigma\cap N  \),
  the relation \( v\leq_{\sigma} v' \) holds if and only if   
    \( \overline{T_{\infty}^X(v)}\supset T_{\infty}^X(v') \).
\end{prop}

\begin{thm}
  Let \( X \) be an affine toric variety and \( x \) a point of \( X \).
  Then the local Nash map:
   \begin{enumerate}
  \item[ ]
   \({{\ell} {\cal N}}:\{ \) local Nash components of \((X,x) \}\)
   \item[ ]
   \ \ \ \ \ \ \ \  \ \ \ \ \ \ \(\to \{ \) essential 
  divisors over \(( X , x) \}\)
  \end{enumerate}
  is bijective.
\end{thm}  

\begin{pf}
  Let \( \sigma \) be the cone defining \( X \).
  We divide the proof into two cases.

{\bf Case 1}: The closure \( \overline{\{x\}} \) is an invariant set.

  In this case, \( \overline{\{x\}} \) is \( orb(\tau) \) for a face 
  \( \tau \) of \( \sigma \) in a 
  neighborhood of \( x \).
   In this neighborhood, \( X=X'\times T' \) for an affine toric variety
  \( X' \) and a torus \( T' \).
  Then, \( \overline{\{x\}}\simeq \{x'\}\times T' \),
  where the point \( x' \) is the closed orbit of \( X' \).
  Therefore, a local Nash component of \( (X,x) \) is of the 
  following type:
\begin{center}
  \( \left( \operatorname{ a\ local\ Nash\ component\ of\ } (X',x')\right) \times {T'}_{\infty} \). 
\end{center}  
  This shows that the number of the local Nash components of \( (X,x) \)
  is that of \( (X',x') \).
  On the other hand, the product of \( T' \)  and 
  a resolution of \( X' \)  
   is a resolution of \( X \) in the neighborhood of \( x \).
  Therefore, an essential divisor over \( (X,x) \) is of  type:
  the product of  \( T' \) and an essential divisor over \( (X',x') \).
  This implies the number of the essential divisors over \( (X,x) \) 
  is
  less than or equal to that over \( (X',x') \).
  Hence we can reduce the problem to the case that \( x \) is the closed 
  orbit.  
    
    Let \( x \) be the closed orbit \( orb (\sigma) \) in \( X \).
  We claim that 
  \[ \pi_{X}^{-1}(x)=\bigcup_{v\in \sigma^o\cap N}
  \overline{T^X_{\infty}(v)}. \]
  For every \( \alpha\in T^X_{\infty}(v) \) with \( v\in \sigma^o\cap
  N \), it follows \( \alpha(0)\in orb(\sigma)=\{x\} \) as we remark 
  in \ref{T}.
  This implies that \( \alpha\in \pi_{X}^{-1}(x) \).
  For  the opposite inclusion,
  it is sufficient to prove that the generic point 
  \( \alpha \) of an irreducible
  component \( C \) of \( \pi_{X}^{-1}(x) \) is contained in 
  \( T^X_{\infty}(v) \) for some \( v\in \sigma^o\cap N \).  
  Let \( \varphi:Y \to X \) be an equivariant resolution.
  Then the induced map \( \varphi_{\infty}:Y_{\infty}\to
  X_{\infty} \)  is surjective (\cite[Proposition 3.2]{i}).
  Therefore, there exists a lifting \( \tilde\alpha\in Y_{\infty} \)
  of \( \alpha \).
  Let \( \tilde\alpha(0)\in orb(\tau) \) for some cone \( \tau \) in 
  the fan of \( Y \).
  Let \( E \) be the closure of \( orb(\tau) \) in \( Y \).
  As \( Y \) is non-singular, \( \pi_{Y}^{-1}(E) \) is irreducible.
  Let \( \beta \) be the generic point of \( \pi_{Y}^{-1}(E) \).
  Since the generic point of \( \pi_{Y}^{-1}(E) \) is  fat (see \ref{fat}),   
  \( \beta(\eta) \) is the generic point of \( Y \)
  therefore  it is in \( T \).
  The inclusion \( \alpha=\varphi_{\infty}(\tilde\alpha)\in
  \varphi_{\infty}(\pi_{Y}^{-1}(E)) \) yields the inclusion 
  \[ C\subset \overline{\varphi_{\infty}(\pi_{Y}^{-1}(E)) } \]
  and this inclusion is an equality,
  because both are irreducible closed subsets  of \( \pi_{X}^{-1}(x) \)
  and \( C \) is an irreducible component of \( \pi_{X}^{-1}(x) \).
  Hence, \( \alpha=\varphi_{\infty}(\beta) \) and therefore \( \alpha(\eta)
  \in T \).
  By this we have a ring homomorphism 
  \[ \alpha^*:\bC[\sigma^{\vee}\cap M]\to K[[t]] \]
  which is extended to 
  \[ \alpha^*:\bC[ M]\to K((t)). \]
  Defining \( v:M\to \bZ \) by \( v(u)=\ord _{t}\alpha^*(x^u) \),
  we obtain \( v\in \sigma\cap N \).
  By \( \alpha(0)=x \), we have \( v\in \sigma^o\cap N \).
  Therefore, we have that 
   \( \alpha\in T_{\infty}^X(v) \).
    
    Now, noting that \( T_{\infty}^{X}(v) \) is not contained in 
    \( (\sing X)_{\infty} \), 
    we see that every irreducible component of \( \pi_{X}^{-1}(x) \) 
    is a local Nash component of \( (X,x) \) and 
    it is a maximal element
    of \( \{\overline{ T_{\infty}^{X}(v)}\mid v\in \sigma^o\cap N\} \).
    Then, by Proposition \ref{dominate}
    the local Nash components of \( (X,x) \) are
    \[ \{ \overline{ T_{\infty}^X (v)}\mid v\ 
    \operatorname{minimal \ in} \sigma^o\cap N\}. \]
    
    On the other hand, 
    an essential divisor over \( (X,x) \) is the same as 
    ``composantes essentielles'' in \cite{B} and the characterization
    theorem of composante essentielle in \cite[\S 2.3]{B} shows that
    \[ \{D_{v}\mid v \ 
    \operatorname{minimal \ in} \sigma^o\cap N\} \]  
    is the set of  composantes essentielles over \( (X,x) \),
    where \( D_{v} \) is the divisor corresponding to the 
    one-dimensional cone \( \bR_{\geq 0}v \).
    (This can be proved also in the similar way as the proof of 
    \cite[Lemma 3.15]{i-k}.)
    This shows the local Nash map is bijective.
    
    Here, we should note that  the proof in \cite[\S 2.3]{B}
    shows that the essential divisors over \( (X,x) \) in the 
    category of all resolutions coincides with that in the category 
    of all equivariant resolutions.

{\bf Case 2}: The closure \( \overline{\{x\}} \) is not an invariant set. 
  
  To prove this case,
  we need the following lemma
\begin{lem}
\label{inclusion}
  Let \( \varphi:Y\to X \) be an equivariant resolution of a toric 
  variety.
	  Let \( orb( \tau )\) be an orbit in \( X \) and \( Z 
  \subset Z' \) irreducible invariant closed subsets of  
  \( \varphi^{-1}(orb( \tau)) \) .
  If \( Z\neq Z' \) then 
  \( Z\cap \varphi^{-1}(\Sigma)\neq Z'\cap \varphi^{-1}(\Sigma)  \) for a 
  subset \( \Sigma\subset orb (\tau) \).  
\end{lem}

\begin{pf}
  As \( Z \) and \( Z' \) are invariant closed subsets of \( 
  \varphi^{-1}(orb( \tau)) \),
  there are lower dimensional toric varieties \( Z_{0}\subset Z'_{0} \)
   such 
  that \( Z\simeq Z_{0}\times orb( \tau )\), \( Z'\simeq Z'_{0}\times 
  orb( \tau) \) and 
  the restrictions of the morphism \( \varphi \) on \( Z,Z' \) are the 
  projections to the second factors.
  Then \( Z\cap \varphi^{-1}(\Sigma)\simeq Z_{0}\times \Sigma \) and 
  \( Z'\cap \varphi^{-1}(\Sigma)\simeq Z'_{0}\times \Sigma  \).
  Hence, \( Z\neq Z' \) implies \( Z_{0}\neq Z'_{0} \) 
  and 
  therefore \( Z\cap \varphi^{-1}(\Sigma)\neq Z'\cap \varphi^{-1}(\Sigma) \).
\end{pf}

  Now we start the proof for Case 2.
  Take the face \( \tau<\sigma \) such that \( x\in orb( \tau) \). 
  Let \( \Sigma= orb( \tau) \cap \overline{\{x\}}\) and let \( x_{\tau} \) 
  be the generic point of \( \overline{orb( \tau)} \).
  Then, we can prove that 
  \[\#\{\mbox{essential divisors over \( (X,x) \)}\}\leq
  \#\{\mbox{essential divisors over \( (X,x_{\tau} )\)}\}.   \]
  In order to prove this, it is sufficient  to prove that for a 
  fixed equivariant resolution \( \varphi:Y\to X \),
  
  \( \#\{\mbox{essential components over \( (X,x) \) on \( Y \)}\} \)
  
  \hskip3truecm
 \( \leq
  \#\{\mbox{essential components over \( (X ,x_{\tau} )\) on \( Y 
  \)}\}\).
  \newline
  Let \( \varphi^{-1}(orb( \tau))=\bigcup_{i=1}^rV_{i} \) be the 
  irreducible decomposition.
  Let \( \Sigma_{i}=V_{i}\cap \varphi^{-1}(\Sigma) \).
  Then, \( \varphi^{-1}(\Sigma)= \bigcup_{i=1}^r\Sigma_{i}\) is the irreducible 
  decomposition. 
  An essential component over \( (X ,x_{\tau} )\) on \( Y 
  \) is one of \( \overline{V_{i}} \), 
  and an essential component over  \( (X,x) \) on \( Y \) is one of
  \( \overline{\Sigma_{i}} \).
  By taking a suitable \( \varphi \) we may assume that 
  \( \overline{V_{i}} \)'s are divisors.
  
 It is sufficient to prove that if \( \overline{V_{i}} \) is not 
  an essential component over \( (X ,x_{\tau} )\) on \( Y 
  \), then \( \overline{\Sigma_{i}} \) is not an essential component over
  \( (X,x) \) on \( Y \).
  If \( \overline{V_{i}} \) is not an essential component over            
  \( (X ,x_{\tau} )\) on \( Y \), there is an equivariant resolution 
  \( \psi:Y'\to X \) such that the center \( \overline{V'_{i}} \) of 
  \( V_{i} \) is strictly contained in an invariant irreducible component \( 
  \overline{V'} \) of \( \overline{\psi^{-1}(x_{\tau})} \).
  Let \( V'_{i}=\overline{V'_{i}}\cap \psi^{-1}(orb( \tau)) \) and 
  \( V'=\overline{V'}\cap \psi^{-1}(orb( \tau)) \).
  Then, by Lemma \ref{inclusion}, the strict inclusion \( V'_{i}\subset
   V' \) 
  yields the strict inclusion 
  \[ (1)\ \ \ \ \ \ \ \ \ \ \ \ \ \  \ \ \ \ \ \ \ \ V'_{i}\cap 
  \psi^{-1}(\Sigma)\subset 
  V'\cap \psi^{-1}(\Sigma)  .\ \ \ \ \ \ \ \ \ \ \ \ \ \ \ \ \ \ \ \ \ \ \ \ \ \ \ \]
  Let   \( g:\tilde Y \to Y \) be 
   an equivariant morphism such that 
  \( \varphi\circ g \) is a resolution of the singularities of \( X \) and 
  there is a morphism \( h: \tilde Y \to Y' \).
  As \( g \) is equivariant and  the minimal invariant closed 
  subset containing \( \overline{\Sigma_{i}} \) is \( \overline{V_{i}} \),
   there is a unique irreducible component \( \widetilde \Sigma_{i} \) of \( 
  g^{-1}(\overline{\Sigma_{i}}) \) mapped onto \( \overline{\Sigma_{i}} \).
  Here, we note that  
  \( \widetilde \Sigma_{i} \subset \overline{V_{i}} \), where we use the 
  same notation for  
  the divisors 
  \( \overline{V_{i}}\subset Y \)  and its proper transform on \( 
  \tilde Y \).  
  Let  \( D \) be an exceptional divisor over \( (X,x) \) whose center on \( Y \)
  is \( \overline{\Sigma_{i}} \).  
  Then the center of \( D \) on \( \tilde Y \) is \( \widetilde \Sigma_{i} \) 
  and therefore the center of \( D \) in \( Y' \) is 
  \( h(\widetilde \Sigma_{i}) \) which is in \( \overline{V'_{i}\cap
  \psi^{-1}(\Sigma)} \).
  By the strict inclusion (1), 
  \( h(\widetilde \Sigma_{i}) \) is contained in another component 
  \( \overline{V'\cap \psi^{-1}(\Sigma)} \).
  Therefore, \( \overline{\Sigma_{i}} \) is not an essential component 
  over \( (X,x) \) on \( Y \).

  Next, we claim that 
  \[ \#\{\mbox{Nash components of \( (X,x) \)}\}=
     \#\{\mbox{Nash components of \( (X,x_{\tau}) \)}\}. \]
  This is proved as follows:
  At a neighborhood of \( x \), \( X\simeq X'\times T' \) and 
  \( orb( \tau)=\{0\}\times T' \), 
  where \( T' \) is a torus of lower dimension, \( X' 
  \) is a suitable toric variety with the closed point orbit \( 0 \).
  We can write \( \Sigma=\{0\}\times \Sigma' \),
  where \( \Sigma'\subset T' \) is an irreducible closed subset.
  Then, \( \pi_{X}^{-1}(orb( \tau))=\pi_{X'}^{-1}(0)\times (T')_{\infty} 
  \) and \(\pi_{X}^{-1}(\Sigma)= \pi_{X'}^{-1}(0)\times \pi_{T'}^{-1}(\Sigma') \).
  Therefore, the numbers of irreducible components of both subsets are 
  the same as the number of irreducible components of 
  \( \pi_{X'}^{-1}(0) \).
  
  Now, using the affirmative answer to the local Nash problem for \( 
  (X,x_{\tau}) \) and the injectivity of the local Nash map,
   we obtain the bijectivity of the local Nash map for \( (X,x) \).

\end{pf}

\section{The local Nash problem for an analytically pretoric singularity}

\noindent
  In \cite{i2}, a pretoric variety is defined and affirmative answer 
  to the Nash problem for a pretoric variety is proved. 
  In this section we introduce an analytically pretoric singularity 
  and give an affirmative answer to the local Nash problem for this 
  singularity.
  A good example of a pretoric variety is a non-normal toric variety,
  while a good example of an analytically pretoric singularity is an
  analytically irreducible 
  quasi-ordinary singularity.
\begin{defn}
\label{pretoric}
  Let \( \o \) be  an integral domain which is the completion of a local ring 
 essentially of finite type over \( k \).
  Let \( X=\spec \o \).
  The closed point of \( X \) is denoted by \( x \).
  A singularity \( (X,x) \) is called an analytically pretoric 
  singularity if the following is satisfied:
  Let \( N=\bZ^n \) and \( M \) the dual of \( N \).
  There exist an \( n \)-dimensional cone  \( \sigma \) 
   in \( N \) and a sublattice \( M'\subset M \) of finite index.
   There is a sequence of injective local homomorphisms
   \[ k[[\sigma^{\vee}\cap M']]\stackrel{\rho^*}\longrightarrow
   \o \stackrel{\nu^*}\longrightarrow k[[\sigma^{\vee}\cap M]], \ \]
  such that
\begin{enumerate}
  \item \( \nu^*\circ\rho^*:k[[\sigma^{\vee}\cap M']]\to  k[[\sigma^{\vee}\cap M]] \ \)
   is the canonical injection,
  \item
  \(  k[[\sigma^{\vee}\cap M]] \) is the integral closure of \( \o \) 
  in its quotient field, and
  \item
   Let \( \nu:\spec k[[\sigma^{\vee}\cap M]]\to \spec \o \) be
   the morphism corresponding to \( \nu^* \).
   The restriction of \( \nu \) onto \( \spec k[[\sigma^{\vee}\cap 
   M]][M] \) is an isomorphism onto its image.
\end{enumerate}  

\end{defn}

\begin{exmp}
  One important example of analytically pretoric singularity is an 
  analytically irreducible quasi-ordinary singularity.
   A quasi-ordinary singularity is first introduced by J. Lipman 
  \cite{lipman1}, \cite{lipman2} and studied by Y-N. Gau \cite{gau}, 
  K. Oh \cite{oh} and P. D. Gonz\'alez P\'erez \cite{perez} and others.

  We call a singularity \( (X,x) \) a {\it quasi-ordinary singularity} if it 
  is a hypersurface singularity in \( (\bC^{n+1},0) \) and 
  there is a finite covering \( \rho: (X, x)\to (\bC^n, 0)  \) whose 
  discriminant locus is contained in a germ wise in a normal crossing divisor on \( \bC^n \).    
  P. D. Gonz\'alez P\'erez \cite{perez} proved that 
  if \( (X,x) \) is an analytically irreducible quasi-ordinary 
  singularity, then it satisfies the conditions of our analytically pretoric singularity.
\end{exmp}

\begin{say}
\label{notation}
  Let \( (X,x) \) be an analytically pretoric singularity. 
  Under the notation in Definition \ref{pretoric},
  we denote \( \spec k[\sigma^{\vee}\cap M] \) and 
  \( \spec k[[\sigma^{\vee}\cap M]] \) by \( W \) and \( 
  {\widehat{W}} \), respectively.
  We denote  \( \spec k[\sigma^{\vee}\cap M'] \) and 
  \( \spec k[[\sigma^{\vee}\cap M']] \) by \( Z \) and \( 
  {\widehat{Z}} \), respectively.
  By the definition of analytically pretoric singularity,
  we obtain the following diagram:
  \[ {\widehat{W}}\stackrel{\nu}\longrightarrow X 
  \stackrel{\rho}\longrightarrow {\widehat{Z}}. \]
  We also obtain that \( \rho\circ\nu \) induces an equivariant 
  morphism \( W\to Z \) of toric varieties.
  Let \( w\in W \)  and \( z\in Z \) be the closed points orbits.
  We denote the closed point of \( \widehat{W} \) and \( \widehat{Z} \)
  by the same symbols \( w \) and \( z \).
  Then they correspond to the point \( x\in X \) by the morphism \( \nu \)
  and \( \rho \).

  As \( T_{\infty}^W(v)\subset \pi_{W}^{-1}(w) \)
  for a point \( v\in \sigma^o\cap N  \), 
  we have \( T_{\infty}^W(v)\subset \pi_{\widehat{W}}^{-1}(w)
  \subset {\widehat{W}} \) by Lemma 
  \ref{fiber}.
  In the same way, for \( v\in \sigma^o\cap N' \), where \( N' \) is 
  the dual of \( M' \),   we obtain that 
  \( T_{\infty}^Z(v)\subset \pi_{\widehat{Z}}^{-1}(z)
  \subset  {\widehat{Z}} \).
\end{say}

\begin{defn}
  For \( v\in \sigma^o\cap N \), define the subset \( T_{\infty}^X(v) \) 
  by the image \( \nu_{\infty}(T_{\infty}^{{W}}(v)) \).
\end{defn}

\begin{lem}
  Let \( (X,x) \) be an analytically pretoric singularity.
  Under the notation in \ref{notation},
  we obtain the following
\begin{enumerate}
  \item[(i)]
  The restriction  
  \( T_{\infty}^W(v) \to T_{\infty}^X(v) \) of \( \nu_{\infty} \)  is bijective for every 
  \( v\in \sigma^o\cap N \).
  \item[(ii)]
  The restriction 
  \( T_{\infty}^W(v) \to T_{\infty}^Z (v) \) 
  of  \( (\rho\circ \nu)_{\infty} \)  is surjective for every 
  \( v\in \sigma^o\cap N \).
\end{enumerate}  
\end{lem}

\begin{pf}
   The surjectivity of (i) follows from the definition of 
  \(  T_{\infty}^X(v)  \).
  The injectivity follows from the valuative criterion of properness, 
   as \( \nu:\widehat{W}\to X \) is proper and the image of \( \eta 
   \) by every arc in  \( T_{\infty}^X(v) \) is in the open set on 
   which \( \nu \) is isomorphic (see (3) in Definition \ref{pretoric}).
  As \( W\to Z \) is the equivariant quotient morphism of toric varieties
  by a finite group \( N'/N \), (ii) follows from \cite[Lemma 5.6, (ii)]{i2}. 
\end{pf}

\begin{lem}
\label{dominant}
  For two points \( v,v'\in \sigma^o\cap N \),
  the following are equivalent
\begin{enumerate}
  \item[(i)]
  \( v\leq_{\sigma}v' \),
  \item[(ii)]
  \( \overline{T_{\infty}^X(v)}\supset T_{\infty}^X(v') \).
\end{enumerate}
\end{lem}

\begin{pf}
  If \(  v\leq_{\sigma}v'  \), then by Proposition \ref{dominate},
   \( \overline{T_{\infty}^W(v)}\supset T_{\infty}^W(v') \).
   Hence, it follows that 
  \[  \overline{\nu_{\infty}(T_{\infty}^W(v))}\supset 
  \nu_{\infty}(T_{\infty}^W(v')), \]
  which implies (ii).
  
  Conversely, if 
  \( \overline{T_{\infty}^X(v)}\supset T_{\infty}^X(v') \),
  then \(  \overline{\rho_{\infty}(T_{\infty}^X(v))}\supset 
  \rho_{\infty}(T_{\infty}^X(v'))  \) which is the inclusion
  \[ \overline{T_{\infty}^Z(v)}\supset T_{\infty}^Z(v'). \]
  Again by Proposition \ref{dominate}, it follows \( v\leq_{\sigma}v' \).
\end{pf}

\begin{lem}
\label{minimal}
  Let \( v \) be a minimal element in \( \sigma^o\cap N \) with 
  respect to the order \( \leq_{\sigma} \).
  Then \( \overline{T_{\infty}^X(v)} \) is a local Nash component of 
  \( (X,x) \).
\end{lem}

\begin{pf}
  As \( T_{\infty}^X(v) \)  is irreducible,
  we can take a Nash component \( C \) of \( (X,x) \) containing 
  \( T_{\infty}^X(v) \) .
  Let \( \alpha \) be the generic point of \( C \),
  then the image \( \alpha(\eta) \) of the generic point \( \eta \) of 
  \( \spec K[[t]] \) is the generic point of \( \widehat{X} \) by Lemma 
  \ref{fat}.  
  Then \( \alpha \) can be uniquely lifted to an arc \( \tilde \alpha:
  \spec K[[t]]\to {\widehat{W}} \) by the valuative criterion of
  properness.
  As  \( \tilde\alpha(\eta) \) is the generic point of \( \widehat{W} \),
  \( \tilde\alpha(\eta) \) is mapped to the generic point of \( W \).
  Then, there exists \( v'\in \sigma^o\cap N \) such that 
  \( \tilde \alpha\in T_{\infty}^{W}(v') \).
  Since \( \overline{\nu_{\infty}T_{\infty}^W(v')}\supset 
  \overline{\{\nu_{\infty}(\tilde\alpha)\}}=\overline{\{\alpha\}}\supset
  T_{\infty}^X(v) \),
  We obtain \[ \overline{T_{\infty}^X(v')}\supset T_{\infty}^X(v). \]
  By Lemma \ref{dominant} and the minimality of \( v \), 
  it follows that \( v=v' \) and 
  \(  C= \overline{T_{\infty}^X(v)} \).     
\end{pf}

\begin{lem}
\label{normalization}
  Let \( \nu:\widehat{W}\to X \) be the normalization of a reduced
  \( k \)-scheme \( X \) and for a singular closed point \( x\in X \), 
  \( \nu^{-1}(x) \) be one closed point \( w \). 
   Then, an essential divisors over \( (X,x) \) is an essential 
   divisors over \( (\widehat{W}, w) \). 
\end{lem}

\begin{pf}
  Let \( E \) be an essential divisor over \( (X,x) \).
  Let \( \psi:\widetilde W \to \widehat{W} \) be a resolution of the 
  singularities of \( \widehat{W} \). 
  Then the composite \( \varphi=\nu\circ \psi:\widetilde W \to X \) is a 
  resolution of the singularities of \( X \) and the center of \( E 
  \) on \( \widetilde W \) is an irreducible component of \( 
  \varphi^{-1}(x)=\psi^{-1}(w) \).  
\end{pf}

\begin{lem}
\label{ess of widehat}
  Let \( w\in W \) be a closed point of a variety and let 
  \( \widehat{W}=\Spec \widehat{\o}_{W,w} \). 
  Denote the closed point of  \( \widehat{W} \) again by \( w \).
  Then an essential divisor over \( (\widehat{W},w) \) is an 
  essential divisor over \( (W,w) \).
\end{lem}

\begin{pf}
  Let \( E \) be an essential divisor over \( (\widehat{W},w) \).
  Then \( E \) is regarded as an exceptional divisor over \( (W,w) \).
  Indeed, for a resolution \( \varphi:Y\to W \) such that \( 
  \varphi^{-1}(w) \) is a divisor,
  the base change \( \varphi':Y\times_{W}\widehat{W}\to \widehat{W} \) 
  is a resolution of the singularities of \( \widehat{W} \) with 
  \( {\varphi'}^{-1}(w)=\varphi^{-1}(w) \).
  As \( E \) appears in \(  {\varphi'}^{-1}(w) \) as a component,
  we can identify \( E \) with the corresponding exceptional divisor over \( (W,w) \).  
  Let \( \psi:Y'\to W \) be any resolution  of the singularities of \( W \)
  and \( \psi':Y'\times_{W}
  \widehat{W}\to \widehat{W} \) the induced resolution of the 
  singularities of \( \widehat{W} \) which is the base change.
 Now, as  \( E \) is an essential divisor over \( (\widehat{W},w) \),
  the center of \( E \) on  \( Y'\times_{W}\widehat{W} \) is an 
   irreducible component of \( {\psi'}^{-1}(w)=\psi^{-1}(w) \). 
\end{pf}

\begin{thm}
\label{main}
  Let \( (X,x) \) be an analytically pretoric singularity.
  Then the local Nash map :
    \begin{enumerate}
  \item[ ]
   \({{\ell} {\cal N}}:\{ \) local Nash components of \((X,x) \}\)
   \item[ ]
   \ \ \ \ \ \ \ \  \ \ \ \ \ \ \( \to\{ \) essential 
  divisors over \(( X , x) \}\)
  \end{enumerate}
  is bijective.

\end{thm}

\begin{pf}
Consider the following diagram:
  $$
\left\{
\begin{array}{c}
\mbox{minimal elements}\\
\mbox{in}\\
\mbox{\( v\in\sigma^o\cap N \)}
\end{array}
\right\}
 \stackrel{\Phi_{1}}\longrightarrow   
    \left\{ 
    \begin{array}{c}
    \mbox{local Nash}\\  
    \mbox{components}\\
    \mbox{ of \( (X,x) \)}
    \end{array} 
    \right\}
 \stackrel{\ell\cal N}\longrightarrow   
  \left\{ 
    \begin{array}{c}
    \mbox{essential}\\  
    \mbox{divisors}\\
    \mbox{ over \( (X,x) \)}
    \end{array} 
    \right\}
$$
$$
\stackrel{\Phi_{2}}\longrightarrow
\left\{ 
    \begin{array}{c}
    \mbox{essential}\\  
    \mbox{divisors}\\
    \mbox{ over \( (\widehat{W},w) \)}
    \end{array} 
    \right\}
 \stackrel{\Phi_{3}}\longrightarrow 
  \left\{ 
    \begin{array}{c}
    \mbox{essential}\\  
    \mbox{divisors}\\
    \mbox{ over \( ({W},w) \)}
    \end{array} 
    \right\}
\stackrel{\Phi_{4}}\longrightarrow     
\left\{
\begin{array}{c}
\mbox{minimal elements}\\
\mbox{in}\\
\mbox{\( v\in\sigma^o\cap N \)}
\end{array}
\right\}.$$

  The map \( \Phi_{1} \) is defined by \( v \mapsto 
  \overline{T_{\infty}^X(v)} \) and it is injective by 
  Lemma \ref{minimal}.
   The local Nash map \( \ell{\cal N} \)  is injective as noted in 
   \ref{locNash}.
  The canonical map \( \Phi_{2} \) is injective by Lemma 
  \ref{normalization}.
  The canonical map \( \Phi_{3} \) is injective by Lemma \ref{ess of widehat}.
  The map \( \Phi_{4} \) sends \( D_{v} \) to \( v \) and it is 
  bijective by Bouvier's characterization of 
  ``composante essentielle''( \cite{B}), where \( D_{v} \) is the 
  invariant divisor \( \overline{orb(\bR_{\geq 0}v)} \).
  Hence the composite of all maps is an injective map from a finite set 
  to itself
   and therefore all 
  maps are bijective.  
\end{pf}

For the final result, we need the following lemma.

\begin{lem}
\label{union}
 Let \( \o \) be  the completion of a local ring 
 essentially of finite type over \( k \) by the maximal ideal.
  Let \( X=\spec \o \).
  The closed point of \( X \) is denoted by \( x \).
Assume that  $X$ is reduced and $X=\bigcup_{i=1}^r X_i$ is the decomposition 
into irreducible components.
If the local Nash map is bijective for $(X_i, x)$ $(i=1,...,r)$,
then the local Nash map is bijective for $(X, x)$.
\end{lem}

\begin{pf}
  Note that \( \pi_{X}^{-1}(x)=\bigcup_{i=1}^r \pi_{X_{i}}^{-1}(x) \).
  First we claim that 
  \[ \{\mbox{local Nash components of } 
  (X,x)\}=\bigsqcup_{i=1}^r\{\mbox{local Nash 
  components of } (X_{i},x)\}. \]
  Let \( C \) be a local Nash component of \( (X,x) \) and 
  \( \alpha \) the generic point of \( C \).
  As \( \alpha(\eta)\in X\setminus \sing X\subset 
  \bigsqcup_{i=1}^r\left( X_{i}\setminus \bigcup_{j\neq i}X_{j}\right), \)
  there is unique \( i \) such that \( \alpha(\eta)\in X_{i} \).
  Then \( \alpha\in \pi_{X_{i}}^{-1}(x) \), 
  therefore \( C\subset \pi_{X_{i}}^{-1}(x)  \) and \( C \) is a 
  local Nash component of \( (X_{i},x) \).
  
  Conversely let \( C \) be a local Nash component of \( (X_{i},x) \) 
  and \( \alpha \) the generic point of \( C \).
  Then \( \alpha(\eta) \) is the generic point of \( X_{i} \) by
  Lemma \ref{fat}.
  Hence, \( \alpha(\eta)\not\in \sing X \).
  Let \( C' \) be a local Nash component of \( (X,x) \) containing \( C \). 
  Then, by the preceding discussion there is unique \( j  \) such 
  that \( C' \) is a local Nash component of
  \( (X_{j},x) \).
  As \( C' \) contains an arc \( \alpha \) satisfying that  \( 
  \alpha(\eta) \) is the generic point of \( X_{i} \), it turns out 
  that \( j=i \).
  Then \( C=C' \) and \( C \) is a local Nash component of \( (X,x) \).
  
  Next we claim that
  \[ \{\mbox{essential divisors over }(X,x)\}\subset \bigsqcup_{i=1}^r\{
  \mbox{essential divisors over } (X_{i},x)\}. \]
  Let \( E \) be an essential divisor over \( (X,x) \), then \( E \) 
  is an exceptional divisor over \( (X_{i},x) \) for some \( i \).
  Let \( \varphi_{i}:Y_{i}\to X_{i} \) be a resolution of the
  singularities of \( X_{i} \).
  Take a resolution \( \varphi_{j}:Y_{j}\to X_{j} \) for each \( 
  j\neq i \).
  Then the composite
  \[ Y:=\bigsqcup_{j=1}^rY_{j}\stackrel{\bigsqcup 
  \varphi_{j}}\longrightarrow \bigsqcup_{j=1}^rX_{j}\to X \]
  is a resolution of the singularities of \( X \).
  As \( E \) is an essential divisor over \( (X,x) \), 
  the center of \( E \) on \( Y \) is an irreducible component of 
  \( \varphi^{-1}(x)=\bigsqcup \varphi_{i}^{-1}(x) \), therefore an 
  irreducible component of \( \varphi_{i}^{-1}(x) \).
  
  Now we obtain the diagram 
    $$
\bigsqcup_{i=1}^i
\left\{
\begin{array}{c}
\mbox{local Nash}\\
\mbox{components}\\
\mbox{of \( (X_{i},x) \)}
\end{array}
\right\}
 =   
    \left\{ 
    \begin{array}{c}
    \mbox{local Nash}\\  
    \mbox{components}\\
    \mbox{ of \( (X,x) \)}
    \end{array} 
    \right\}
 \stackrel{\ell\cal N}\longrightarrow   
  \left\{ 
    \begin{array}{c}
    \mbox{essential}\\  
    \mbox{divisors}\\
    \mbox{ over \( (X,x) \)}
    \end{array} 
    \right\}
$$

\[\subset 
\bigsqcup_{i=1}^i
\left\{
\begin{array}{c}
\mbox{essential}\\
\mbox{divisors}\\
\mbox{over \( (X_{i},x)  \)}
\end{array}
\right\}.
  \]
  
  Since the local Nash components of $(X_i, x)$ correspond bijectively to the
  essential components over $(X_i,x)$ for each $i$,
  the composite of all injections of above diagram is bijective.
    Therefore all maps are bijective.
\end{pf}

Now we obtain the following final statement.

\begin{cor}
  Let \( (X,x) \) be a quasi-ordinary singularity. 
  Then the local Nash map for \( (X,x) \) is bijective.
\end{cor}
 
\begin{pf}
  A quasi-ordinary singularity \( (X,x) \) is decomposed into 
  analytically irreducible quasi-ordinary singularities 
  \( (X_{i}, x) \) \( (i=1,..,r) \).
As each \( (X_{i},x) \)  is an analytically pretoric singularity,
the statement follows from Theorem \ref{main} and Lemma \ref{union}.
\end{pf}


\makeatletter \renewcommand{\@biblabel}[1]{\hfill#1.}\makeatother

\end{document}